\documentstyle[12pt]{article} \vsize=28.7truecm \hsize=23truecm
  \def\vbar{\mathchoice{\vrule height2.3ptdepth-.3ptwidth.12pt\kern-
 .10pt}
    {\vrule height6.3ptdepth-.3ptwidth.11pt\kern-.11pt}
    {\vrule height5.1ptdepth-.30ptwidth.8pt\kern-.8pt}
    {\vrule height4.1ptdepth-.24ptwidth.6pt\kern-.7pt}}
 \def\fudge{\mathchoice{}{}{\mkern.5mu}{\mkern.8mu}}
 \def\bbc#1#2{{\rm \mkern#2mu\vbar\mkern-#2mu#1}}
 \def\bbb#1{{\rm I\mkern-2.5mu #1}}
 \def\bba#1#2{{\rm #1\mkern-#2mu\fudge #1}}
 \def\bb#1{{\count4=`#1 \advance\count4by-64 \ifcase\count4\or\bba
 A{12.5}\or
    \bbb B\or\bbc C{5}\or\bbb D\or\bbb E\or\bbb F \or\bbc G{5}\or\bbb
 H\or
    \bbb I\or\bbc J{3}\or\bbb K\or\bbb L \or\bbb M\or\bbb N\or\bbc
 O{5} \or
    \bbb P\or\bbc Q{5}\or\bbb R\or\bbc S{4.2}\or\bba T{10.5}\or\bbc
 U{5}\or
    \bba V{12}\or\bba W{16.5}\or\bba X{11}\or\bba Y{11.7}\or\bbay
 Z{7.5}\fi}}

 \def \n {\noindent}

 \def \R {\bb R} % reals
\begin{document}
 \vspace{0.5cm}

\title{{\bf
 A characterization of the Riesz distribution}\\
 {\small (Running title: \textbf{A characterization of the Riesz distribution})}}
\author {A. Hassairi\footnote{Corresponding author.
 \textit{E-mail address: Abdelhamid.Hassairi@fss.rnu.tn}}$\;$  , S. Lajmi and R.
 Zine
\\{\footnotesize{\it
Facult\'e des Sciences, Universit\'e de Sfax, B.P.802, Sfax, Tunisie.}}\\
    {\footnotesize{\it  }}}

 \date{}
 \maketitle

\n $\overline{\hspace{16.5cm}}$\vskip0.3cm

\n {\small {\bf Abstract}}\\

{\small Bobecka and Wesolowski (2002) have shown that, in the
Olkin and Rubin characterization of the Wishart distribution (See
Casalis and Letac (1996)), when we use the division algorithm
defined by the quadratic representation and replace the property
of invariance by the existence of twice differentiable densities,
we still have a characterization of the Wishart distribution. In
the present work, we show that, when we use the division algorithm
defined by the Cholesky decomposition, we get a characterization
of the Riesz
distribution.\\

\n {\small {\it Keywords:} Symmetric cone, division algorithm,
Wishart distribution, Riesz distribution, Beta-Riesz
distribution, functional equation}\\
$\overline{\hspace{16.5cm}}$\vskip1cm

\section{Introduction}
A remarkable characterization of the gamma distribution, due to
Luckacs (1955), says that if $U\;$and $V\;$are two independent non
Dirac and non negative random variables such that $U+V\;$is a.s.
positive, then $\frac{U}{U+V}\;$and $U+V\;$are independent if and
only if $U\;$and $V\;$have the gamma distribution with the same
scale parameter. The classical multivariate version of this
characterization concerns the Wishart distribution on the cone of
symmetric positive matrices. In this case, there is not a single
way to define the quotient of two matrices. For instance if $Y$ is
a positive definite matrix, one can for example, use the quadratic
representation, that is write
$Y=Y^{\frac{1}{2}}Y^{\frac{1}{2}}\;$and define the ratio $X\;$by
$Y$ as $Y^{-\frac{1}{2}}XY^{-\frac{1}{2}},\;$or use the Cholesky
decomposition $Y=TT^{\ast }$, where $T\;$is a lower triangular
matrix and define the ratio$\;$as ($T^{-1})X(T^{-1})^{\ast }.\;$A
general definition of a division algorithm will be given in
Section 2, however these two examples are the most usual and most
important. In 1962, Olkin and Rubin have shown that, if $U\;$and
$V\;$are two independent random variables valued in the cone of
symmetric non negative matrices such that $U+V\;$is a.s. positive
definite, then, independently of the choice of the division
algorithm, the quotient of $U\;$by $U+V\;$is independent of
$U+V\;$and its distribution is invariant by the orthogonal group
if and only if $U\;$and $V\;$have the Wishart distribution. This
result has been extended to the Wishart distribution on any
symmetric cone by Casalis and Letac (1996). Recently, Bobecka and
Wesolowski (2002) have given another characterization of the
Wishart distribution without any invariance assumption for the
quotient. More precisely, they have shown that if we use the
division algorithm defined by the quadratic representation and we
replace in the Olkin and Rubin theorem the condition of invariance
for the distribution of the quotient by the existence of twice
differentiable densities, then we still have a characterization of
the Wishart distribution. The present paper gives a parallel
result which starts from the observation that, when the condition
of invariance of the distribution of the ratio by the orthogonal
group is dropped, the characterization in the Bobecka and
Wesolowski way is not independent of the choice of the division
algorithm. We show that, when we use the division algorithm
defined by the Cholesky decomposition, we get a characterization
of the Riesz distribution introduced by Hassairi and Lajmi (2001).
Our method of proof is based on some functional equations
depending on the triangular group. These equations are more
involved then the ones used in the characterization of the Wishart
distribution, their solutions are expressed in terms of the
generalized power. Our results will be presented in the framework
of the Riesz distribution on the symmetric cone of a simple
Euclidean Jordan algebra. This will enable us to use some
technical results established in the book of Faraut-Kor\'{a}nyi
(1994) and in Hassairi et al. (2001, 2005). However, to make the
paper accessible to a reader who is not familiar with the theory
of Jordan algebras, we will give a particular emphasis to the cone
of positive definite symmetric matrices.

\section{Riesz distributions}

We first review some facts concerning Jordan algebras and their
symmetric cones. Our notations are the ones used in the book of
Faraut-Kor\'{a}nyi (1994). Let us recall that a Euclidean Jordan
algebra is a Euclidean space $E$ with scalar product $<x,y>$ and a
bilinear map
$$E\times E\longrightarrow E,\ (x,y)\longmapsto xy$$
\noindent called Jordan product such that, for all $x,\ y,\ z$ in
$E$,

\noindent i) $xy=yx$,

\noindent ii) $<x,yz>=<xy,z>$,

\noindent iii) there exists $e$ in $E$ such that $ex=x$,

\noindent iv) $x(x^2y)=x^2(xy)$, where we used the
abbreviation $x^2=xx$.\bigskip

An Euclidean Jordan algebra is said to be simple if it does not
contain a nontrivial ideal. Actually to each Euclidean simple
Jordan algebra, one attaches the set of Jordan squares
$$\overline{\Omega}=\{x^2;\ x\in E\}.$$
\noindent Its interior $\Omega$ is a symmetric cone, i.e., a cone
which is

\noindent i) self dual, i.e., $\Omega=\{x \in E;\ <x,y>\
>0\ \forall y \in \overline{\Omega}\setminus \{0\} \}$.

\noindent ii) homogeneous, i.e., the subgroup $G(\Omega)$ of the
linear group $GL(E)$ of linear automorphisms which preserves
$\Omega$ acts transitively on $\Omega$.

\noindent iii) salient, i.e., $\Omega$ does not contain a line.
Furthermore, it is irreducible in the sense that it is not the
product of two cones.\bigskip

Let now $x$ be in $E$. If $L(x)$ is the endomorphism of $E$;
$y\longmapsto xy$ and $ P(x)=2L(x)^2-L(x^2)$, then $L(x)$ and
$P(x)$ are symmetric for the Euclidean structure of $E$ and the
map $x\longmapsto P(x)$ is called the quadratic representation of
$E$.\bigskip

An element $c$ of $E$ is said to be idempotent if $c^2=c$, it is a
primitive idempotent if furthermore $c\neq 0$ and is not the sum
$t+u$ of two non null idempotents $t$ and $u$ such that
$t.u=0$.\bigskip

A Jordan frame is a set $\{c_1,...,c_r\}$ of primitive idempotents
such that $\displaystyle \sum_{i=1}^rc_i=e$ and
$c_ic_j=\delta_{ij}c_i$, for $1\leq i,j\leq r$. It is an
important result that the size $r$ of such a frame is a constant
called the rank of $E$.\bigskip

If $c$ is a primitive idempotent of $E$, the only possible
eigenvalues of $L(c)$ are $0,\frac{1}{2}$ and $1$. The
corresponding eigenspaces are respectively denoted by $E(c,0),\
E(c,\frac{1}{2})$ and $E(c,1)$ and the decomposition
$$E=E(c,0)\oplus E(c,\frac{1}{2})\oplus E(c,1)$$
\noindent is called the Peirce decomposition of $E$ with respect
to $c$.

Suppose now that $(c_i)_{1\leq i \leq r}$ is a Jordan frame in $E$
and let, for $1\leq i,j \leq r$,
$$ E_{ij}=\left\{
\begin{array}{l}
E(c_i,1)=\R c_i\ \ \ \ \ \ \ \ \ \ \textrm{if}\ i=j  \\
E(c_i,\frac{1}{2})\cap E(c_j,\frac{1}{2})\ \ \ \ \  \textrm{if}\
i\neq j.
\end{array}
\right. $$

\noindent Then (See Faraut-Kor\'{a}nyi (1994), Theorem IV.2.1) we
have $E=\displaystyle\oplus_{i\leq j}E_{ij}$ and the dimension of
$E_{ij}$ is, for $i\neq j$, a constant $d$ called the Jordan
constant. It is related to the dimension $n$ and the rank $r$ of
$E$ by the relation $n=r+r(r-1)\frac{d}{2}.$\bigskip

For $1\leq k\leq r$, let $P_k$ denote the orthogonal projection on
the Jordan subalgebra
$$E^{(k)}=E(c_1+...+c_k,1),$$
\noindent $\det^{(k)}$ the determinant in the subalgebra $E^{(k)}$
and, for $x$ in $E$, $\Delta_k(x)=\det^{(k)}(P_k(x))$. Then
$\Delta_k$ is called the principal minor of order $k$ with respect
to the Jordan frame $(c_i)_{1\leq i\leq r}$. For
$s=(s_1,...,s_r)\in\R^r$, and $x$ in $\Omega$, we write
$$\Delta_s(x)=\Delta_1(x)^{s_1-s_2}\Delta_2(x)^{s_2-s_3}......\Delta_r(x)^{s_r}.$$This
is the generalized power function. Note that, if $x=\displaystyle
\sum_{i=1}^r\lambda_ic_i$, then
$\Delta_s(x)=\lambda_1^{s_1}\lambda_2^{s_2}...\lambda_r^{s_r}$ and
that $\Delta_s(x)=(\det x)^p$ if $s=(p,...,p)$ with $p \in \R$. It
is also easy to see that
$\Delta_{s+s'}(x)=\Delta_s(x)\Delta_{s'}(x)$. In particular, if $m
\in \R \ $ and $s+m=(s_1+m,...,s_r+m)$, we have
$\Delta_{s+m}(x)=\Delta_{s}(x)(\det x)^m$. \bigskip

As we have mentioned above, one may suppose that $E$ is the
algebra of real symmetric matrices with rank $r$. In this case,
the Jordan product $xy$ of two symmetric matrices $x$ and $y$ is
defined by $\frac{1}{2}(x.y+y.x)$ where $x.y$ is the ordinary product
of the matrices $x$ and $y$, the cone $\Omega$ is the cone of
positive definite matrices, $\overline{\Omega}$ is the cone of
symmetric non negative matrices and $d=1$. If $x= (x_{ij} )_{1\leq
i,j \leq r}$ is an $(r,r)-$symmetric positive definite matrix, and
if, for $1\leq k\leq r$, we denote $P_k(x)=(x_{ij} )_{1\leq i,j
\leq k}$ and $\Delta_k (x)=\det(x_{ij})_{1\leq i,j \leq k}$, the
generalized power is the function on $\Omega$ defined by
$$\Delta_s(x)=\Delta_1(x)^{s_1-s_2}\Delta_2(x)^{s_2-s_3}......\Delta_r(x)^{s_r}.$$

The definition of the Riesz distribution on a symmetric cone
$\Omega$ relies on the notion of generalized power. In fact, for
$\sigma$ is in $\Omega$ and $s$ such that, for all $i$,
$s_{i}>(i-1)\frac{d}{2}$, the measure on $\Omega$

$$
R(s,\sigma )(dx)=\frac{1}{\Gamma _{\Omega }(s)\Delta _{s}(\sigma
^{-1})}e^{-<\sigma ,x>}\Delta _{s-\frac{n}{r}}(x)
{\mathbf{1}}_{\Omega }(x)dx
$$
\n where $\Gamma_{\Omega }(s)=(2\pi
)^{\frac{n-r}{2}}\displaystyle\prod_{j=1}^r\Gamma
(s_{j}-(j-1)\frac{d}{2})$, is a probability distribution. It is
called the Riesz distribution with parameters $s$ and $\sigma$.

We come now to the general definition of a division algorithm in a symmetric cone $%
\Omega $. Let $G$ be the connected component of the identity in
$G(\Omega )$.\ A division algorithm is defined as a measurable map
$g$ from $\Omega $ into $G$ such that, for all $y\;$in $\Omega
,\;$ $g(y)(y)=e$.\bigskip

As in the case of symmetric matrices, we will introduce two
important division algorithms, the first is based on the quadratic
representation $ x\mapsto P(x^{-\frac{1}{2}} )$ and the second
algorithm takes its values in the triangular group $T$. For the
definition of $T$, we need to introduce some other facts
concerning a Jordan algebra. For $x$ and $y$ in $E$, let $x\Box y$
denote the endomorphism of $E$ defined by
\begin{equation}\label{d1}
  x\Box y=L(xy)+[L(x),L(y)]=L(xy)+L(x)L(y)-L(y)L(x).
\end{equation}

If $c$ is an idempotent and if $z$ is an element of
$E(c,\frac{1}{2})$,
$$\tau_c(z)=\exp(2z\Box c)$$is called a Frobenius transformation, it is an element of the group $G$.\bigskip

Given a Jordan frame $(c_i)_{1\leq i \leq r}$, the subgroup of $G$

$$
T=\left\{ \tau _{c_{1}}(z^{(1)}).....\tau
_{c_{r-1}}(z^{(r-1)})P(\displaystyle\sum_{i=1}^r a_ic_i) ,\
a_i>0,\ z^{(j)}\in \displaystyle\bigoplus_{k=j+1}^r E_{jk}\right\}
$$
is called the triangular group corresponding to the Jordan frame
$(c_i)_{1\leq i \leq r}$. It is an important result
(Faraut-Kor\'{a}nyi, p.113, Prop VI.3.8) that the symmetric cone
$\Omega$ of the algebra $E$ is parameterized by the set

\begin{equation}\label{a1}
  E_+=\{u=\displaystyle\sum_{i=1}^ru_ic_i+\sum_{i<j}u_{ij},\ u_i>0\}
\end{equation}
More precisely, if

\begin{equation}\label{a2}
 t_u=\tau_{c_{1}}(z^{(1)}).....\tau_{c_{r-1}}(z^{(r-1)})P(\displaystyle\sum_{i=1}^r u_ic_i)
\end{equation}
where $z_{ij}=\frac{u_{ij}}{u_i},\ i<j$ and
$z^{(j)}=\displaystyle\sum_{k=j+1}^rz_{jk}$, then the map
$u\longmapsto t_u(e)$ is a bijection from $E_+$ into $\Omega$ with
a Jacobian equal to $2^r\displaystyle\prod_{i=1}^r
u_i^{1+d(r-i)}$. Also, for all $x$ in $E$, we have
$$\Delta_k(t_u(x))=u_1^2...u_k^2\Delta_k(x)=\Delta_k(t_u(e))\Delta_k(x).$$
It is shown that, for each $b$ in $\Omega $, there exists a unique
$t$ in the triangular group $T$ such that $b=t(e)$. Hence the map
\begin{equation}\label{h1}
 g:\Omega \longrightarrow T;\ b\longmapsto t^{-1}
\end{equation}

\noindent realizes a division algorithm.
\section{Characterization of the Riesz distribution}

In this section we state and prove our main result which may be
seen as an extension of the result by Bobecka and Wesolowski
(2002) concerning the ordinary Wishart distribution on symmetric
matrices. More precisely these authors use the division algorithm
defined by the quadratic representation to characterize the
Wishart distribution, we use the division algorithm defined by the
triangular group to characterize the Riesz distribution.
\bigskip

\n{\bf Theorem 3.1.} {\it Let $X$ and $Y$ be two independent Riesz
random variables $X\sim R(s,\sigma)$ and

\n $Y\sim R(s^\prime,\sigma)$. If we set $V=X+Y$ and
$U=g(X+Y)(X)$, then

\noindent i) $V$ is a Riesz random variable $V\sim
R(s+s^\prime,\sigma)$ and is independent of $U$

\noindent ii) The density of $U$ with respect to the Lebesgue
measure is
$$
\frac{1}{B_{\Omega }(s,s^{\prime
})}\Delta_{s-\frac{n}{r}}(x)\Delta _{s^{\prime
}-\frac{n}{r}}(e-x){\mathbf{1}}_{\Omega \cap (e-\Omega )}(x),
$$where $B_{\Omega }(s,s^{\prime })$ is the beta function defined
on the symmetric cone $\Omega $ (See Faraut-Kor\'{a}nyi, 1994,
p.130) by
$$
B_{\Omega }(s,s^{\prime })=\frac{\Gamma _{\Omega }(s)\Gamma
_{\Omega }(s^{\prime })}{\Gamma _{\Omega }(s+s^{\prime })}.$$}

\n \textbf{Proof.}  Consider the transformation:
$\Omega\times\Omega\longrightarrow (\Omega\cap
(e-\Omega))\times\Omega;\ (x,y)\longmapsto(u,v)$, where
$v=x+y=te,\ t\in T$ and $u=t^{-1}(x)$. Its Jacobian is $\det
t^{-1}=(\det t(e))^{-\frac{n}{r}}=(\det v)^{-\frac{n}{r}}$.

\bigskip
\n The density of probability of $(U,V)$ with respect to the
Lebesgue measure is then given by
$$ \frac{1}{\Gamma_\Omega(s) \Gamma_\Omega (s^{ \prime}) \Delta_{s+s^{ \prime}}(\sigma^{-1})}\  \Delta_{s-\frac{n}{r}}(t(u))
\Delta_{s^{\prime}-\frac{n}{r}}(t(e-u))\ (\det v)^{\frac{^n}{r}}\
e^{- < \sigma,v >} {\mathbf{1}}_{K}(u,v)$$where $K$ is defined by
$$K=\{(u,v)/u\in \Omega\cap(e-\Omega)\ \textrm{and}\ v\in
\Omega\}.$$ Using equality $(3.5)$ in Hassairi and Lajmi $(2001)$,
this density may be written as
$$\frac{\Gamma_\Omega (s+s^{ \prime})}{\Gamma_\Omega (s)\Gamma_\Omega (s^{ \prime})}
\ \Delta_{s-\frac{n}{r}}(u) \
\Delta_{s^{\prime}-\frac{n}{r}}(e-u)\frac{1}{\Gamma_\Omega (s+s^{
\prime})\Delta_{s+s^{ \prime}}(\sigma^{-1})} \ e^{- < \sigma,v >}\
\Delta_{s+s^{ \prime}-\frac{n}{r}}(v){\mathbf{1}}_{K}(u,v).$$

\n From this we deduce that $U$ and $V$ are independent and that
$V\sim R(s+s^\prime,\sigma)$. Furthermore the distribution of $U$
is concentrated on $\Omega\cap (e-\Omega)$ with a density equal to
$$\frac{\Gamma_\Omega (s+s^{ \prime})}{\Gamma_\Omega (s)\Gamma_\Omega (s^{ \prime})}
\Delta_{s-\frac{n}{r}}(u)
\Delta_{s^{\prime}-\frac{n}{r}}(e-u){\mathbf{1}}_{\Omega\cap(e-\Omega)}(u).$$
$$\eqno\Box$$

Note that the distribution of the random variable $U=g(X+Y)(X)$ is
called beta-Riesz distribution with parameters $s$ and $s^\prime$
(See Hassairi et al (2005)). In the case of symmetric matrices,
the random variable $U$ is nothing
but$$U=(T^{-1})X(T^{-1})^\ast,$$ where $T$ is a lower triangular
matrix with positive diagonal such that$$X+Y=TT^\ast.$$

\n{\bf Theorem 3.2.} {\it Let $b\longmapsto g(b)$ be the division
algorithm defined by (\ref{h1}). Let $X$ and $Y$ be independent
random variables valued in $\Omega$ with strictly positive twice
differentiable densities. Set $V=X+Y$ and $U=g(V)(X)$. If $U$ and
$V$ are independent then there exist $s,\ s^\prime\in \R^r $; $s_i>(i-1)\frac{d}{2}$, $\ s^\prime_i>(i-1)\frac{d}{2}\ $ for all $i$, and
$\sigma\in\Omega$ such that $X\sim R(s,\sigma)$ and $Y\sim
R(s^{\prime},\sigma)$}.\bigskip

The proof of this theorem relies on the resolution of two
functional equations given in the following theorems which are
interesting in their own rights. The proofs of these theorems are
given in Section 4.\\

\n {\bf Theorem 3.3.}  {\it Let
$a:\Omega\cap(e-\Omega)\longrightarrow \R$ and
$g:\Omega\longrightarrow\R$ be functions such that, for any $x\in
\Omega\cap(e-\Omega)$ and $t\in T$,
\begin{equation}\label{b1}
 a(x)=g(tx)-g(t(e-x)).
\end{equation}

\noindent Assume that $g$ is differentiable, then there exist
$p\in \R^r$ and $c\in \R$ such that, for any $x\in
\Omega\cap(e-\Omega)$ and $y\in \Omega,$
$$a(x)=\log\Delta_p(x)-\log\Delta_p(e-x),\
g(y)=\log\Delta_p(y)+c.$$}

\noindent {\bf Theorem 3.4.}  {\it Let
$a_1:\Omega\cap(e-\Omega)\longrightarrow \R$ and $a_2,\
g:\Omega\longrightarrow \R$ be functions satisfying

\begin{equation}\label{b8}
  a_1(x)+a_2(te)=g(tx)+g(t(e-x)),
\end{equation}

\n for any $x\in \Omega\cap (e-\Omega)$ and $t\in T$. Assume that
$g$ is twice differentiable then there exist $p^{\prime}\in \R^r$,
 $\delta\in E$ and $c_1,\ c_2,\ c_3\in \R$ such
that for any $x\in \Omega\cap (e-\Omega)$ and $y\in \Omega$,
$$g(y)=\log\Delta_{p^\prime}(y)+<\delta,y>+c_1$$
$$a_1(x)=\log\Delta_{p^\prime}(x)+\log\Delta_{p^\prime}(e-x)+c_2,$$
$$a_2(y)=2\log\Delta_{p^\prime}(y)+<\delta,y>+c_3,$$
where $2c_1=c_2+c_3.$}\bigskip

\n{\bf Proof of Theorem 3.2.} We again use the transformation:
$\Omega\times\Omega\longrightarrow (\Omega\cap
(e-\Omega))\times\Omega;\ (x,y)\longmapsto(u,v)$, where
$v=x+y=te,\ t\in T$ and $u=t^{-1}(x)$. Let $f_X,\ f_Y,\ f_U$ and
$f_V$ be the densities of $X,\ Y,\ U$ and $V$, respectively. Then,
since $(X,Y)$ and $(U,V)$ have independent components, we have
that, for all $u\in \Omega\cap (e-\Omega)$ and $v\in \Omega$,

\begin{equation}\label{z1}
  f_U(u)f_V(v)=(\det
v)^{\frac{n}{r}}f_X(tu)f_Y(t(e-u))
\end{equation}

\n Taking logarithms in (\ref{z1}) we get

\begin{equation}\label{z2}
  g_1(u)+g_2(te)=g_3(tu)+g_4(t(e-u)),
\end{equation}

\noindent where

\begin{equation}\label{z3}
  g_1(u)=\log f_U(u),
\end{equation}

\begin{equation}\label{z4}
  g_2(v)=\log f_V(v)-\frac{n}{r}\log\det v,
\end{equation}

\begin{equation}\label{z5}
 g_3(x)=\log f_X(x),
\end{equation}

\n and

\begin{equation}\label{z6}
  g_4(y)=\log f_Y(y).
\end{equation}

\noindent Inserting $e-u$ for $u$ in (\ref{z2}) gives

\begin{equation}\label{z7}
  g_1(e-u)+g_2(te)=g_3(t(e-u))+g_4(tu).
\end{equation}

\n Subtracting (\ref{z7}) from (\ref{z2}), we obtain

\begin{equation}\label{z8}
  g_1(u)-g_1(e-u)=[g_3(tu)-g_4(tu)]-[g_3(t(e-u))-g_4(t(e-u))].
\end{equation}

\noindent Define$$a(u)=g_1(u)-g_1(e-u),\
g=g_3-g_4,$$then,$$a(u)=g(t(u))-g(t(e-u)).$$Now according to
Theorem 3.3, we obtain $$a(u)=\log\Delta_p(u)-\log\Delta_p(e-u),\
g(v)=\log\Delta_p(v)+c,$$for $p=(p_1,...,p_r)\in \R^r$ and $c\in
\R$.

\n Hence

\begin{equation}\label{z9}
  g_3(v)=g_4(v)+g(v)=g_4(v)+\log\Delta_p(v)+c.
\end{equation}

\noindent Inserting (\ref{z9}) back into (\ref{z2}) gives
$$g_1(u)+g_2(te)=\log\Delta_p(u)+\log\Delta_p(te)+c+g_4(t(u))+g_4(t(e-u)),$$which
can be rewritten in the form

\begin{equation}\label{z10}
  a_1(u)+a_2(te)=g_4(t(u))+g_4(t(e-u)),
\end{equation}

\noindent where

\begin{equation}\label{z11}
  a_1(u)=g_1(u)-\log\Delta_p(u),
\end{equation}

\begin{equation}\label{z12}
 a_2(te)=g_2(te)-\log\Delta_p(te)-c.
\end{equation}

\noindent Hence, by Theorem 3.4, it follows that
$$g_4(v)=\log\Delta_{p^\prime}(v)+<\delta,v>+c_1,$$
$$\ \ \ \ \ \ \ \ a_1(u)=\log\Delta_{p^\prime}(u)+\log\Delta_{p^\prime}(e-u)+c_2,$$
$$\ \ a_2(v)=2\log\Delta_{p^\prime}(v)+<\delta,v>+c_3,$$
for some $\delta\in E$, ${p^\prime}\in \R^r$ and $c_1,\ c_2,\
c_3\in \R$ such that $c_3=2c_1-c_2.$\bigskip

\n This with (\ref{z6}) imply$$\log
f_Y(y)=\log\Delta_{p^\prime}(y)+<\delta,y>+c_1, $$that
is$$f_Y(y)=e^{c_1}\ e^{<\delta,y>}\Delta_{p^\prime}(y).$$Since
$f_Y$ is a probability  density, it follows that
$\sigma=-\delta\in \Omega$ and ${s^\prime}=
{p^\prime}+\frac{n}{r}$ is such that
${s^\prime}_i>(i-1)\frac{d}{2},\ \forall 1\leq i\leq r$ so that
$Y\sim R({s^\prime},\sigma)$.\bigskip

\n Now by (\ref{z9}), we get

\begin{eqnarray*} g_3(v) & = &
\log\Delta_p(v)+c+\log\Delta_{p^\prime}(v)+<\delta,v>+c_1
\\ & = &
\log\Delta_{p+{p^\prime}}(v)+<\delta,v>+c+c_1.
\end{eqnarray*}

\noindent From (\ref{z5}) it follows that

\begin{eqnarray*} f_X(x) & = &
\Delta_{p+{p^\prime}}(x)e^{<\delta,x>}e^{c+c_1} \\ & = &
\Delta_{p+{s^\prime}-\frac{n}{r}}(x) e^{-<\sigma,x>}e^{c+c_1}
\end{eqnarray*}

\noindent which implies that $p_i+{s^\prime}_i>(i-1)\frac{d}{2},\
\forall 1\leq i\leq r$ and consequently $X\sim R(s,\sigma)$ where
$s=p+{s^\prime}$.

\section{Proofs of Theorems 3.3 and 3.4}
\subsection{Proof of Theorem 3.3}

For the proof of Theorem 3.3, we need to establish some
preliminary results.\bigskip

\n {\bf Proposition 4.1. }{\it i) The map
$\varphi:\Omega\longrightarrow\R;\ x\longmapsto \log\Delta_k(x)$
is differentiable on $\Omega$ and
$$\nabla\log\Delta_k(x)=(P_k(x))^{-1}.$$

\n ii) The differential of the map $x\longmapsto (P_k(x))^{-1}$ is
$-P((P_k(x))^{-1}).$

}\bigskip

\n{\bf Proof.} We first observe that if $c$ is an idempotent of
$E$ and $x$ is in $E(c,1)$, then $$P(x) P(c)= P(x)$$ \n In fact,
let $h= h_1+h_{12}+h_0 $ be the Peirce decomposition of an element
$h$ in $E$ with respect to $c$. As $x\in E(c,1)$, we have that
$$ P(x)(h_{12})= P(x)(h_{0})=0,$$
\n and it follows that
$$ P(x)(h)= P(x)(h_{1})= P(x)P(c)(h)$$

\n i) Let $\Omega_k=P_k(\Omega)$ and consider the map
$\psi:\Omega_k\longrightarrow\R;\ x\longmapsto\log\det^{(k)}(x)$.
Then $\varphi(x)=\psi\circ P_k(x)$. Since $\nabla\log\det
x=x^{-1}$, we have

\begin{eqnarray*}
  (\log\Delta_k(x))^{\prime}(h) & = &{\psi}^{\prime}(P_k(x))(P_k(x))^\prime(h),\ \forall h\in E\\ & = &
  <(P_k(x))^{-1},P_k(h)>\\ &=&  <(P_k(x))^{-1},h>
 \end{eqnarray*}

 \n ii) We have that the differential in $x$ of the map $\beta:x\longmapsto x^{-1}$ is
$-P(x^{-1})$. As $(P_k(x))^{-1}=\beta\circ P_k(x),\
 \forall x\in \Omega$, then the differential in $ x\in\Omega$ of the map $x\mapsto (P_k (x))^{-1}$
 is equal to $ - P((P_k (x))^{-1})\circ P_k $ which is also equal to $ - P((P_k (x))^{-1})$\bigskip

 \n {\bf Proposition 4.2.} {\it Let $u$ be in $E_+$ and let $x=t_u(e).$ Then, for $1\leq k\leq r,$
 $$(P_k(x))^{-1}=t_u^{{\ast}^{-1}}(\sum_{i=1}^kc_i).$$}

\n{\bf Proof.} We know, from Faraut-Kor\'{a}nyi (Proposition
VI.3.10), that for $x$ in $E$ and $t_u$ in $T$,
 $$P_k(t_u(x))=t_{P_k(u)}(P_k(x)).$$Using (\ref{a2}) we can write

 \begin{eqnarray*}
  P_k(x) & = &\tau({\widetilde{z}}^{(1)})...\tau({\widetilde{z}}^{(k-1)})P(\sum_{i=1}^ku_ic_i)(c_1+...+c_k)\\ & = &
  \tau({\widetilde{z}}^{(1)})...\tau({\widetilde{z}}^{(k-1)})(\sum_{i=1}^ku_i^2c_i)
 \end{eqnarray*}

\n where ${\widetilde{z}}^{(j)}=\displaystyle\sum_{l=j+1}^kz_{jl}$
and $z_{jl}=\frac{u_{jl}}{u_j},\ j<l.$

\n Hence

\begin{equation}\label{a4}
  (P_k(x))^{-1}=\tau(-{\widetilde{z}}^{(1)})^\ast...\tau(-{\widetilde{z}}^{(k-1)})^\ast\ (\sum_{i=1}^k\frac{1}{u_i^2}c_i)
\end{equation}

\n On the other hand
$$t_u^{{\ast}^{-1}}(\sum_{i=1}^kc_i)=\tau(-z^{(1)})^\ast...\tau(-z^{(r-1)})^\ast\
(\sum_{i=1}^k\frac{1}{u_i^2}c_i).$$

\n Let us recall that if $c$ is an idempotent, $z$ is in
$E(c,\frac{1}{2})$ and $x_1,\ x_{12},\ x_0$ are the Peirce
components of $x$ with respect to $c$. Then the Peirce components
of $y=\tau_c(z)^\ast(x)$ are

\begin{equation}\label{a5}
  \left\{ \begin{array}{l} y_1=2c[z(zx_0)+zx_{12}]+x_1,  \\
y_{12}=2zx_0+x_{12}, \\ y_0=x_0. \end{array} \right.
\end{equation}

\n This, after some elementary calculations, implies that, for
$k+1\leq j\leq r-1$, $z^{(j)}\in E(c_j,\frac{1}{2})$ and
$\displaystyle\sum_{i=1}^k\frac{1}{u_i^2}c_i\in E(c_j,0)$
$$\tau(-z^{(j)})^\ast \
(\sum_{i=1}^k\frac{1}{u_i^2}c_i)=\sum_{i=1}^k\frac{1}{u_i^2}c_i,\
\forall k+1\leq j\leq r-1, $$ and it follows that,
$$t_u^{{\ast}^{-1}}(\sum_{i=1}^kc_i)=\tau(-z^{(1)})^\ast...\tau(-z^{(k)})^\ast\
(\sum_{i=1}^k\frac{1}{u_i^2}c_i).$$Using again (\ref{a5}), we have
that
$$t_u^{{\ast}^{-1}}(\sum_{i=1}^kc_i)=\tau(-z^{(1)})^\ast...\tau(-z^{(k-1)})^\ast\
(\sum_{i=1}^k\frac{1}{u_i^2}c_i).$$From this, we deduce that, for
$1\leq j\leq k-1$ and $a\in
E(c_1+...+c_k,1)$,$$\tau(-z^{(j)})^\ast
a=\tau(-\widetilde{z}^{(j)})^\ast a.$$Finally, we conclude that

\begin{eqnarray*}
  t_u^{{\ast}^{-1}}(\sum_{i=1}^kc_i) & = &\tau(-z^{(1)})^\ast...\tau(-z^{(k-1)})^\ast\
(\sum_{i=1}^k\frac{1}{u_i^2}c_i)\\ & = &
  \tau(-\widetilde{z}^{(1)})^\ast...\tau(-\widetilde{z}^{(k-1)})^\ast\
(\sum_{i=1}^k\frac{1}{u_i^2}c_i)\\ & = & (P_k(x))^{-1}.
 \end{eqnarray*}

\bigskip

For the proof of Theorem 3.3, we also need the following result
for which a proof is given in Hassairi and Lajmi (2001). For
$h=\displaystyle\sum_{i=1}^rh_ic_i+\sum_{i<j}h_{ij}$ in $E$, we
denote $ \overline{h}=\displaystyle\sum_{i=1}^rh_ic_i$ and
$h^{(j)}=\displaystyle\sum_{k=j+1}^rh_{jk},\ 1\leq j\leq
r-1.$\bigskip

\n {\bf Lemma 4.3. }{\it i) For all $1\leq j\leq r-1$, the map
from $E_+$ into $L(E)$; $u\longmapsto \tau(u^{(j)})$ is
differentiable on $E_+$. Its differential in a point $u$ is the
map $$h\longmapsto K_j^h(u)=2h^{(j)}\Box c_j+4(h^{(j)}\Box
c_j)(u^{(j)}\Box c_j)$$

\n ii) The map $H:E_+\longrightarrow T;\ u\longmapsto t_u$ is
differentiable on $E_ +$ and

\begin{eqnarray*}
  H^\prime(u)(h) & = &\ \displaystyle\sum_{j=1}^{r-1}\tau(u^{(0)})\tau(u^{(1)})...
\tau(u^{(j-1)})K_j^h(u)\tau(u^{(j+1)})...\tau(u^{(r-1)})P(\overline{u})\\
& + &
2\tau(u^{(1)})...\tau(u^{(r-1)})[2L(\overline{u})L(\overline{h})-L(\overline{u}\overline{h})]
\end{eqnarray*}

\n where $\tau(u^{(0)})=Id.$}\bigskip

\n Note that for $u=e$, we get

\begin{equation}\label{a6}
  K_j^h(e)=2h^{(j)}\Box c_j\ \textrm{and}\ H^\prime(e)(h)
  =\displaystyle\sum_{j=1}^{r-1}K_j^h(e)+2L(\overline{h})=
   2[\sum_{j=1}^{r-1}h^{(j)}\Box
c_j+L(\overline{h})]
\end{equation}

We are now in position to prove Theorem 3.3.\bigskip

\n{\bf Proof of Theorem 3.3.} Differentiating (\ref{b1}) with
respect to $x$ gives

\begin{equation}\label{b2}
  a^{\prime}(x)=t^\ast g^\prime(tx)+t^\ast g^\prime(t(e-x)).
\end{equation}

\n Setting $x=\frac{e}{2}$ and replacing $t$ by $2t$ in (\ref{b2})
give
\begin{equation}\label{b3}
  g^\prime(te)=t^{{\ast}^{-1}}b,\ \textrm{where}\ b=\frac{1}{4}a^\prime(\frac{e}{2})=g^\prime(e).
\end{equation}

\n Taking $t=Id$ in (\ref{b2}) gives
$$a^\prime(x)=g^\prime(x)+g^\prime(e-x).$$Inserting this identity
back into (\ref{b2}) gives
$$t^{{\ast}^{-1}}g^\prime(x)-g^\prime(tx)=-[t^{{\ast}^{-1}}g^\prime(e-x)-g^\prime(t(e-x))].$$
For any $s\in (0,1)$ and $x=te$, we
have$$g^\prime(sx)=g^\prime(ste)=\frac{1}{s}t^{{\ast}^{-1}}b=\frac{1}{s}g^\prime(x).$$
Hence for any $s\in
(0,1)$,$$t^{{\ast}^{-1}}g^\prime(x)-g^\prime(tx)=-s[t^{{\ast}^{-1}}g^\prime(e-sx)-g^\prime(t(e-sx))].$$
Consequently, on letting $s\longrightarrow 0$, we obtain for any
$x\in \Omega\cap (e-\Omega)$ and $t\in T$, that
\begin{equation}\label{b4}
  g^\prime(tx)=t^{{\ast}^{-1}}g^\prime(x).
\end{equation}

\n Now consider $u\in E_+$. Then (\ref{b1}) can be rewritten as

\begin{equation}\label{b5}
  a(x)=g(t_u(x))-g(t_u(e-x)).
\end{equation}

\n  To differentiate (\ref{b5}) with respect to $u$, let us
consider the functions
$$H: E_+\longrightarrow L(E); \ u\longmapsto t_u\ \textrm{and}\ \alpha: L(E)\longrightarrow E;\ f\longmapsto f(x),
$$

\n so that for $u\in E_+$, we can write
$g(t_u(x))=(g\circ\alpha\circ H )(u).$

\n From Lemma 4.3, we have that

\begin{eqnarray*} (g\circ\alpha\circ H )^\prime(u)(h) & = &
g^\prime(t_u(x))\alpha^\prime(H(u))H^\prime(u)(h),\forall h\in E\\
& = & g^\prime(t_u(x))\alpha(H^\prime(u)(h))\\ & = &
<g^\prime(t_u(x)),H^\prime(u)(h)(x)>.
 \end{eqnarray*}

\n Then, for all $h\in E$,
$$<g^\prime(t_u(x)),H^\prime(u)(h)(x)>=<g^\prime(t_u(e-x)),H^\prime(u)(h)(e-x)>.$$
Hence by (\ref{b4}) we
get$$<t_u^{{\ast}^{-1}}g^\prime(x),H^\prime(u)(h)(x)>=<t_u^{{\ast}^{-1}}g^\prime(e-x),H^\prime(u)(h)(e-x)>,$$
and for any $s\in (0,1)$,

\begin{eqnarray*}
<t_u^{{\ast}^{-1}}g^\prime(sx),H^\prime(u)(h)(sx)> & = &
<t_u^{{\ast}^{-1}}g^\prime(x),H^\prime(u)(h)(x)>\\ & = &
<t_u^{{\ast}^{-1}}g^\prime(e-sx),H^\prime(u)(h)(e-sx)>.
 \end{eqnarray*}

\n Letting $s\longrightarrow 0$, we get
$$<t_u^{{\ast}^{-1}}g^\prime(x),H^\prime(u)(h)(x)>=
<t_u^{{\ast}^{-1}}g^\prime(e),H^\prime(u)(h)(e)>.$$Note that if
$u=e$, we have for all $x\in \Omega\cap (e-\Omega)$,

\begin{equation}\label{b6}
 <g^\prime(x),H^\prime(e)(h)(x)>= <b,H^\prime(e)(h)(e)>,\ \forall
h\in E.
\end{equation}

Let $b=\displaystyle\sum_{i=1}^rb_ic_i+\sum_{i<j}b_{ij}$ be the
Peirce decomposition of $b$ defined in (\ref{b3}) with respect to
the Jordan frame $(c_i)_{1\leq i\leq r}$. We will show that
$b=\displaystyle\sum_{i=1}^rb_ic_i$. In order to do so, we
consider $x=\frac{e}{2}+\varepsilon x_{ij}$ with $\varepsilon\in
]0,\frac{1}{2}[$ and $x_{ij}\in E_{ij}$ such that $\|x_{ij}\|=1$.
It is easy to see that $x=t_v(e)$, where
$v=\displaystyle\sum_{k=1,k\neq
j}^r\frac{c_k}{\sqrt{2}}+\sqrt{\frac{1}{2}-\varepsilon^2}\
c_j+\sqrt{2}\varepsilon x_{ij} \in E_+$. This implies that $x\in
\Omega$ and $e-x=\frac{e}{2}-\varepsilon x_{ij}$ is also in
$\Omega$.\bigskip

\n Inserting $h=h_{ij}$ in (\ref{a6}), we obtain
$$H^\prime(e)(h)=2h_{ij}\Box c_i.$$With the notations used above
and from (\ref{d1}), we have
$$H^\prime(e)(h)(e)=h_{ij},$$and
$$H^\prime(e)(h)(x)=\frac{1}{2}h_{ij}+\varepsilon<h_{ij},x_{ij}>c_j.$$
It follows by (\ref{b3}) that

\begin{eqnarray*} <g^\prime(x),H^\prime(e)(h)(x)> & = &
<t_v^{{\ast}^{-1}}b,H^\prime(e)(h)(x)>\\ & = &
<b,t_v^{-1}H^\prime(e)(h)(x)>. \end{eqnarray*}

\n After some standard calculation using (\ref{a2}), we get
$$t_v^{-1}H^\prime(e)(h)(x)=\frac{1}{\sqrt{1-2\varepsilon^2}}h_{ij}.$$Hence
according to (\ref{b6}), we have for all $h_{ij}\in E_{ij}$
$$<b_{ij},h_{ij}>=\frac{1}{\sqrt{1-2\varepsilon^2}}<b_{ij},h_{ij}>.$$
From this we readily deduce that $b_{ij}=0\ \textrm{for all}\
1\leq i<j\leq r$.\bigskip

\n Now let $x=t(e)$, with $t\in T$ be in $\Omega$. Then from
(\ref{b3}), we can write

\begin{eqnarray*} g^\prime(x) & = & t^{{\ast}^{-1}}b\\ & = &
t^{{\ast}^{-1}}(\displaystyle\sum_{i=1}^rb_ic_i)\\ & = &
t^{{\ast}^{-1}}(\displaystyle\sum_{k=1}^r(b_k-b_{k+1})(c_1+...+c_k))
\end{eqnarray*}

\n where $b_{r+1}=0$.\bigskip

\n And using Proposition 4.2, we get
$$g^\prime(x)=\displaystyle\sum_{k=1}^r(b_k-b_{k+1})(P_k(x))^{-1}.$$By
the statement (i) of Proposition 4.1, we have

$$g(x)=
\displaystyle\sum_{k=1}^r(b_k-b_{k+1})\log\Delta_k(x)+c,\ c\in
\R$$

\begin{equation}\label{b7}
  =\log\Delta_p(x)+c\ \ \ \ \ \ \ \ \ \ \ \ \ \ \ \ \ \ \ \ \ \
\end{equation}

\n where $p=(b_1,...,b_r)$.\bigskip

\n Finally, inserting (\ref{b7}) into (\ref{b1}) we
get$$a(x)=\log\Delta_p(x)-\log\Delta_p(e-x).$$
\subsection{Proof of Theorem 3.4}
For the resolution of the functional equation (\ref{b8}), we will
use second derivatives and the following result.

\n {\bf Proposition 4.4.} {\it Let
$x=\displaystyle\sum_{i=1}^rx_ic_i+\sum_{i<j}x_{ij}$ be the Peirce
decomposition with respect to $(c_i)$ of an element $x$ of $E$.
Then, for $q=(q_1,...,q_r)$ in $\R^r,$
$$\displaystyle\sum_{i=1}^rq_ix_ic_i+\sum_{i<j}q_jx_{ij}=[q_rP(c_1+...+c_r)+
\displaystyle\sum_{k=1}^{r-1}(q_k-q_{k+1})P(c_1+...+c_k)](x).$$}

\n{\bf Proof.} We use induction on the rank $r$ of the algebra.
The result is obvious for $r=2$, in fact we have

\begin{eqnarray*} q_1x_1c_1+q_2x_2c_2+q_2x_{12} & = & q_2(x_1c_1+x_2c_2+x_{12})+(q_1-q_2)x_1c_1\\ & = &
[q_2P(c_1+c_2)+(q_1-q_2)P(c_1)](x).
\end{eqnarray*}

\noindent Suppose the result true for a Jordan algebra with a rank
equal to $r-1$. We easily verify that
$$\displaystyle\sum_{i=1}^rq_ix_ic_i+\sum_{i<j}q_jx_{ij}=
q_r(\displaystyle\sum_{i=1}^{r}x_ic_i+\sum_{i<j}x_{ij})+
\displaystyle\sum_{i=1}^{r-1}(q_i-q_r)x_ic_i+\sum_{1\leq i<j\leq
r-1}(q_j-q_r)x_{ij}.$$ As, for
$a=\displaystyle\sum_{i=1}^{r}a_ic_i$, we have
$P(a)x=\displaystyle\sum_{i\leq j}a_ia_j x_{ij}$, then
$$\displaystyle\sum_{i=1}^rq_ix_ic_i+\sum_{i<j}q_jx_{ij}= q_rP(c_1+...+c_r)(x) +$$
$$[(q_{r-1}-q_r)
P(c_1+...+c_{r-1})+\displaystyle\sum_{k=1}^{r-2}(q_k-q_{k+1})P(c_1+...+c_k)]
(\displaystyle\sum_{i=1}^{r-1} x_ic_i+\sum_{1\leq i<j\leq
r-1}x_{ij})$$This with some standard calculation give the
result.\bigskip

\n {\bf Proof of Theorem 3.4.} Differentiating (\ref{b8}) with
respect to $x$ gives

\begin{equation}\label{b9}
  a_1^\prime(x)=t^\ast g^\prime(tx)-t^\ast g^\prime(t(e-x)).
\end{equation}

\n Differentiating (\ref{b9}) once again with respect to $x$ gives

\begin{equation}\label{b10}
  a_1^{\prime\prime}(x)=t^\ast g^{\prime\prime}(tx)t+t^\ast
g^{\prime\prime}(t(e-x))t. \end{equation}

\n Substitute $x=\frac{1}{2}e$ in the above equation then replace
$t$ by $2t$ to get

\begin{equation}\label{b11}
 g^{\prime\prime}(te)=t^{{\ast}^{-1}}Bt^{-1}
\end{equation}

\n where
$B=\frac{1}{8}a^{\prime\prime}(\frac{e}{2})=g^{\prime\prime}(e).$\bigskip

\n Observe that taking $t=Id$ in (\ref{b10}) gives

\begin{equation}\label{b12}
a_1^{\prime\prime}(x)=g^{\prime\prime}(x)+g^{\prime\prime}(e-x).
\end{equation}

\n Inserting this identity back into (\ref{b10}), we get
$$g^{\prime\prime}(x)-t^\ast
g^{\prime\prime}(tx)t=-[g^{\prime\prime}(e-x)-t^\ast
g^{\prime\prime}(t(e-x))t].$$For any $s\in (0,1)$, change $x$ to
$sx$ in the above equation, and use the fact that by (\ref{b11}),
we have$$g^{\prime\prime}(sx)=s^{-2}g^{\prime\prime}(x),$$then on
multiplication by $s^2$, we deduce that for any $s\in (0,1)$
$$g^{\prime\prime}(x)-t^\ast
g^{\prime\prime}(tx)t=-s^2[g^{\prime\prime}(e-sx)-t^\ast
g^{\prime\prime}(t(e-sx))t].$$Letting $s\longrightarrow 0$, we
obtain, for any $x\in \Omega\cap (e-\Omega)$ and $t\in T$

\begin{equation}\label{b13}
  g^{\prime\prime}(tx)=t^{{\ast}^{-1}}g^{\prime\prime}(x) t^{-1}.
\end{equation}

Let us now consider $u\in E_+$, then (\ref{b9}) can be rewritten
in the form

\begin{equation}\label{b14}
  a_1^{\prime}(x)=t_u^\ast g^{\prime}(t_u(x))-t_u^\ast
g^{\prime}(t_u(e-x)). \end{equation}

\n  In order to differentiate (\ref{b14}) with respect to $u$, we
introduce the functions$$H:E_+\longrightarrow L(E);\ u\longmapsto
t_u,\ \pi:L(E)\longrightarrow L(E);\ f\longmapsto f^\ast\
\textrm{and}\ \phi:L(E)\longrightarrow E;\ f\longmapsto f(x).$$

\n We have $$(\pi\circ H)^\prime(u)(h)=H^\prime(u)(h)^\ast,\
\forall h\in E.$$

\n As, for any $u\in E_+$, we have
$g^\prime(t_u(x))=(g^\prime\circ \phi\circ H)(u)$, then
$$(g^\prime\circ \phi\circ
H)^\prime(u)(h)=g^{\prime\prime}(t_u(x))H^\prime(u)(h)(x).$$We
also consider the functions$$\Psi:E_+\longrightarrow L(E)\times
E;\ u\longmapsto(t_u^\ast,g^\prime(t_u(x))),\ \textrm{where}\ x\
\textrm{is fixed},\ \textrm{and}\ \zeta:L(E)\times
E\longrightarrow E;\ (f,z)\longmapsto f(z). $$

\n Then one can easily see that
$$\Psi^\prime(u)(h)=\left(H^\prime(u)(h)^\ast,g^{\prime\prime}(t_u(x))H^\prime(u)(h)(x)\right),$$
and it follows that

\begin{eqnarray*} (\zeta\circ \Psi)^\prime(u)(h) & = &
\zeta^\prime(\Psi(u))\Psi^\prime(u)(h)\\ & = &
\zeta^\prime(t_u^\ast,g^\prime(t_u(x)))(H^\prime(u)(h)^\ast,
g^{\prime\prime}(t_u(x))H^\prime(u)(h)(x))\\ & = &
\zeta(t_u^\ast,g^{\prime\prime}(t_u(x))H^\prime(u)(h)(x))+\zeta(H^\prime(u)(h)^\ast,g^\prime(t_u(x)))\\
& = & t_u^\ast
g^{\prime\prime}(t_u(x))H^\prime(u)(h)(x)+H^\prime(u)(h)^\ast(g^\prime(t_u(x))).
 \end{eqnarray*}

\n Now differentiating (\ref{b14}) with respect to $u$ gives
$$t_u^\ast
g^{\prime\prime}(t_u(x))H^\prime(u)(h)(x)+H^\prime(u)(h)^\ast(g^\prime(t_u(x)))=
t_u^\ast
g^{\prime\prime}(t_u(e-x))H^\prime(u)(h)(e-x)+H^\prime(u)(h)^\ast(g^\prime(t_u(e-x))).$$
If we make $u=e$, we get

\begin{equation}\label{b15}
  g^{\prime\prime}(x)H^\prime(e)(h)(x)+H^\prime(e)(h)^\ast
g^\prime(x)=
g^{\prime\prime}(e-x)H^\prime(e)(h)(e-x)+H^\prime(e)(h)^\ast
g^\prime(e-x). \end{equation}

\noindent That is
$$(g^{\prime\prime}(x)+g^{\prime\prime}(e-x))H^\prime(e)(h)(x)+
H^\prime(e)(h)^\ast (g^\prime(x)-g^\prime(e-x))=
g^{\prime\prime}(e-x))H^\prime(e)(h)(e).$$ Using (\ref{b9}) and
(\ref{b12}), this becomes
$$a_1^{\prime\prime}(x)H^\prime(e)(h)(x)+H^\prime(e)(h)^\ast
a_1^\prime(x)=g^{\prime\prime}(e-x))H^\prime(e)(h)(e).$$On the
other hand, we know that

$$\left \{ \begin{array}{l}
a_1^{\prime\prime}(e-x)=a_1^{\prime\prime}(x)\\
a_1^{\prime}(e-x)=-a_1^{\prime}(x). \end{array} \right.$$Hence
$$a_1^{\prime\prime}(x)H^\prime(e)(h)(e-x)-H^\prime(e)(h)^\ast
a_1^\prime(x)=g^{\prime\prime}(x)H^\prime(e)(h)(e).$$Given
(\ref{b12}), we obtain $$H^\prime(e)(h)^\ast
a_1^\prime(x)=g^{\prime\prime}(e-x)H^\prime(e)(h)(e-x)-g^{\prime\prime}(x)H^\prime(e)(h)(x),\
\forall h\in E. $$Setting $h=e$ in this equality, we obtain by
(\ref{a6}),

 \begin{equation}\label{b16}
  a_1^\prime(x)=g^{\prime\prime}(e-x)(e-x)-g^{\prime\prime}(x)(x).
\end{equation}

\noindent Therefore, using (\ref{b16}) and the fact that, for all
$x$ in $\Omega\cap (e-\Omega)$,
$a_1^\prime(x)=g^\prime(x)-g^\prime(e-x)$, the equality
(\ref{b15}) becomes
\begin{equation}\label{b17}
  g^{\prime\prime}(x)H^\prime(e)(h)(x)-H^\prime(e)(h)^\ast g^{\prime\prime}(x)(x)=
  g^{\prime\prime}(e-x)H^\prime(e)(h)(e-x)-H^\prime(e)(h)^\ast g^{\prime\prime}(e-x)(e-x).
\end{equation}

\noindent For $s\in (0,1)$ change $x$ to $sx$, to obtain $$
g^{\prime\prime}(x)H^\prime(e)(h)(x)-H^\prime(e)(h)^\ast
g^{\prime\prime}(x)(x)=s[
  g^{\prime\prime}(e-sx)H^\prime(e)(h)(e-sx)-H^\prime(e)(h)^\ast g^{\prime\prime}(e-sx)(e-sx)].$$
Letting $s\longrightarrow 0$, implies that for any $x\in
\Omega\cap (e-\Omega)$ and $h\in E$,

\begin{equation}\label{b18}
 H^\prime(e)(h)^\ast g^{\prime\prime}(x)(x)=g^{\prime\prime}(x)H^\prime(e)(h)(x).
\end{equation}

\noindent If $x=\frac{e}{2}$, then, using the fact that
$g^{\prime\prime}(\frac{e}{2})=4B$, we can write

\begin{equation}\label{b19}
  H^\prime(e)(h)^\ast B(e)=B H^\prime(e)(h)(e),\ \forall h\in E.
\end{equation}

\noindent Now let
$$B(e)=\displaystyle\sum_{i=1}^r\lambda_ic_i+\sum_{i<j}\lambda_{ij}$$be
the Peirce decomposition of $B(e)$ with respect to the Jordan
frame $(c_i)_{1\leq i \leq r}$. We will show that
$B(e)=\displaystyle\sum_{i=1}^r\lambda_ic_i$.\bigskip

\n In fact, substituting $h=\displaystyle\sum_{i=1}^rh_ic_i$ in
(\ref{b19}) gives
$$(\displaystyle\sum_{i=1}^rh_ic_i)B(e)=\displaystyle\sum_{i=1}^rh_i
B(c_i)$$In particular, if $h=c_i$, we get$$B(c_i)=c_iB(e),\
\forall 1\leq i \leq r.$$Consequently

\begin{equation}\label{b20}
  B(c_i)=\lambda_ic_i+\frac{1}{2}[\displaystyle\sum_{k=1}^{i-1}\lambda_{ki}+
\displaystyle\sum_{k=i+1}^r\lambda_{ik}]
\end{equation}

\noindent Using (\ref{b19}) again for $h=h_{ij}$, we get

\begin{eqnarray*} B(h_{ij}) & = & H^\prime(e)(h)^\ast B(e)\\ & =
& 2(h_{ij}\Box c_i)^\ast B(e)\\ & = &  2(c_i\Box h_{ij})B(e)\\ & =
& <\lambda_{ij},h_{ij}> c_i+ \lambda_j
h_{ij}+2h_{ij}[\displaystyle\sum_{k=1}^{j-1}\lambda_{kj}+
\displaystyle\sum_{k=j+1}^r\lambda_{jk}],\ \forall i<j.
 \end{eqnarray*}

\noindent By symmetry of $B$, it follows
that$$<B(h_{ij}),c_j>=\frac{1}{2}<h_{ij},\lambda_{ij}>=0,\ \forall
i<j.$$ This implies that $\lambda_{ij}=0,\ \forall i<j.$\bigskip

\noindent Hence, we have

$$\left \{ \begin{array}{l}
B(c_i)=\lambda_ic_i,\ \forall 1\leq i\leq r.\\ B(h_{ij})=\lambda_j
h_{ij},\
\forall 1\leq i<j\leq r.\\
B(e)=\displaystyle\sum_{i=1}^r\lambda_ic_i.\end{array}
\right.$$Let $x=
\displaystyle\sum_{i=1}^rx_ic_i+\displaystyle\sum_{i<j}x_{ij}$ be
the Peirce decomposition of $x$. By Proposition 4.4, we get

\begin{eqnarray*} B(x) & = &
\displaystyle\sum_{i=1}^r\lambda_ix_ic_i+\displaystyle\sum_{i<j}\lambda_jx_{ij}
\\ & = &
[\lambda_rP(c_1+...+c_r)+\displaystyle\sum_{k=1}^{r-1}(\lambda_k-\lambda_{k+1})P(c_1+...+c_k)](x),
 \end{eqnarray*}

\noindent that is
$$B=\lambda_rP(c_1+...+c_r)+\displaystyle\sum_{k=1}^{r-1}(\lambda_k-\lambda_{k+1})P(c_1+...+c_k).$$

As for each $x$ in $\Omega$, there exists a unique $t$ in the
triangular group $T$ such that $x=t(e)$, then using (\ref{b11})
and the fact that, for all $x$ in $E$ and for all $g\in G$,
$P(gx)=gP(x)g^\ast,$ we can write

\begin{eqnarray*} g^{\prime\prime}(x) & = &
t^{{\ast}^{-1}}Bt^{-1} \\ & = &
\lambda_rt^{{\ast}^{-1}}P(c_1+...+c_r)t^{-1}+\displaystyle\sum_{k=1}^{r-1}
(\lambda_k-\lambda_{k+1})t^{{\ast}^{-1}}P(c_1+...+c_k)t^{-1}\\ & =
&
\lambda_rP(t^{{\ast}^{-1}}(c_1+...+c_r))+\displaystyle\sum_{k=1}^{r-1}(\lambda_k-\lambda_{k+1})P(t^{{\ast}^{-1}}(c_1+...+c_k))
\end{eqnarray*}

\noindent Using Proposition 4.2, we can write
$$g^{\prime\prime}(x)=\lambda_rP(x^{-1})+\displaystyle\sum_{k=1}^{r-1}(\lambda_k-\lambda_{k+1})
P((P_k(x))^{-1}).$$This, invoking Proposition 4.1, implies that
$$g^{\prime}(x)=-\lambda_r
x^{-1}-\displaystyle\sum_{k=1}^{r-1}(\lambda_k-\lambda_{k+1})(P_k(x))^{-1}+\delta$$where
$\delta\in E$.

\noindent And by Proposition 4.1, we have that

$$g(x)=-\lambda_r\log\det
x-\displaystyle\sum_{k=1}^{r-1}(\lambda_k-\lambda_{k+1})\log\Delta_k(x)+<\delta,x>+c_1$$where
$c_1\in \R$.

\noindent Hence
$$g(x)=\log\Delta_{p^{\prime}}(x)+<\delta,x>+c_1$$where
$p^{\prime}=(p^{\prime}_1,...,p^{\prime}_r)=(-\lambda_1,...,-\lambda_r)$.

\noindent As $a_1^\prime(x)=g^\prime(x)-g^\prime(e-x)$ (See
(\ref{b9})), we obtain $$a_1(x)=
\log\Delta_{p^\prime}(x)+\log\Delta_{p^\prime}(e-x)+c_2$$where
$c_2\in \R$.

\noindent Finally, from (\ref{b8}) we get
$$a_2(te)=2\log\Delta_{p^\prime}(te)+2c_1-c_2+<\delta,te>.$$For
any $y\in \Omega$ there exists a unique $t\in T$ such that $y=te$.
Then
$$a_2(y)=2\log\Delta_{p^\prime}(y)+<\delta,y>+c_3,$$where
$c_3=2c_1-c_2.$\\

\noindent $\mathbf{Acknowledgements}$ We are very grateful to
Jacek Wesolowski for discussions on the subject of this paper.
 \end{document}